\documentclass[12pt,reqno]{amsart}
  
 \usepackage{comment}
\usepackage{tikz}
\usepackage{mathrsfs}
\usepackage{latexsym}
\usepackage{graphicx}
\usepackage{amssymb}
 \graphicspath{{Figures/}}
 \usepackage[show]{ed}

\renewcommand{\Im}{{\operatorname{Im}\,}}

\renewcommand{\epsilon}{\varepsilon}

\newcommand{\R}{{\mathbb R}}

\newcommand{\supp}{{\operatorname{supp\,}}}
\renewcommand{\phi}{\varphi}

\newtheorem{maintheo}{Theorem}
\newtheorem{cor}{Corollary}[section]

\newtheorem{lem}[cor]{Lemma}

\newtheorem{prop}[cor]{Proposition}
\numberwithin{equation}{section}
\newtheorem{rem}{Remark}

\usepackage{hyperref}

\title[Improved Lower Bound]
{Improved Lower Bound for Analytic Schr\"odinger Eigenfunctions in Forbidden Regions}

\author{Xianchao Wu}
\address{Department of Mathematics, Wuhan University of Technology, Wuhan, Hubei, China}
\email{xianchao.wu@whut.edu.cn} 

\date{}

\begin{document}

\maketitle
\begin{abstract}
The point of this paper is to improve the reverse Agmon estimate discussed in \cite{TW} with assuming that the Schrodinger operator $P(h) = - h^2 \Delta_g + V - E(h)$, $E(h)\to E$ as $h\to 0^+$, is analytic on a compact, real-analytic Riemannian manifold $(M,g)$. In this paper, by considering a Neumann problem with applying Poisson representation and exterior mass estimates on hypersurfaces, we can prove an improved reverse Agmon estimate on a hypersurface.

\end{abstract}

\section{Introduction}

Let $(\mathcal{M}, g)$ be a compact, real-analytic $n$-dimentional Riemannian manifold and $V(x) \in C^{\omega}(\mathcal{M};\R)$ be a real-analytic potential. Assume that $E$ is a regular value of $V$ so that $dV |_{V=E} \neq 0$.
The corresponding classically {\em forbidden region} is denoted as
\begin{equation} \label{allowable}
\Omega_E := \{ x \in \mathcal{M}; V(x) > E \},
\end{equation}
and the {\em allowed region} is the complement $\Omega_E^c = \{x \in \mathcal{M}; V(x) \leq E \}$
with boundary $C^{\omega}$ hypersurface (ie. boundary caustic)
\begin{equation}\label{caustic}
\Lambda_E := \{ x \in \mathcal{M};  V(x) = E \}. \end{equation}

Consider the Schr\"odinger equation
\begin{equation}\label{Schrodinger}
P(h) u_h=0, 
\end{equation}
where $P(h):= - h^2 \Delta_g + V(x) - E(h)$ and $\{u_h\}$ are $L^2$-normalized eigenfunctions with eigenvalues $E(h) \to E$ as $h\to 0^+$. The {\em Agmon metric} is given by
$$ g_E (x) := (V(x)-E)_{+} \, g(x).$$
The degenerate metric $g_E$ is supported in the forbidden region $\Omega_E$ and we denote the corresponding Riemannian distance function by $d_E: \Omega_E \times \Omega_E \to \R^+.$  By a slight abuse of notation, we define the associated distance function to $\Lambda_E$ by

\begin{equation} \label{agmondistance}
d_E(x):= d_{E}(x, \Lambda_E) = \inf_{y \in \Lambda_E} d_{E}(x,y), \quad x \in \Omega_E.
\end{equation}

In \cite{TW} ($(\mathcal{M},g,V)$ is only required to be smooth), under the control and monotonicity assumptions (see \cite[Definitions 1  and 2]{TW}), by applying Carleman estimate to pass across the caustic hypersurface \cite[Theorem 1]{TW} authors prove that for any $\epsilon >0$ and $h \in (0,h_0(\epsilon)],$

\begin{equation} \label{partial}
\| e^{\tau_0 d_E/h} \, u_h \|_{H_h^1 (A(\delta_1,\, \delta_2))} \geq C(\epsilon,\delta_1, \delta_2) \, e^{-\beta(\epsilon)/h}, 
\end{equation}
where $ \beta(\epsilon) = O(\epsilon)$
as $\epsilon \to 0^+$, $A(\delta_1,\delta_2)\subset \Omega_E$ is an annular domain near the boundary (precise definitions refer to \cite{TW}) and constant $\tau_0\geq 1$.

As authors pointed out in \cite{TW}, \eqref{partial} is a partial reverse Agmon estimate, our objective in this paper is to get an improved result in the case where $(\mathcal{M}, g, V)$ is analytic. This is precisely the point of inequality \eqref{result}. We note that the assumption that $(\mathcal{M}, g, V)$ is $C^\omega$ is necessary in this article since real-analyticity allows for accuracy up to exponential errors in $h$ in the pseudodifferential calculus, whereas in the $C^\infty$ case, one can only work to $O(h^\infty)$-error. Since the eigenfunctions decay exponentially in $h$ in forbidden regions, the usual $C^\infty$ semiclassical calculus of operators is not accurate enough to deal with these functions in a rigorous fashion. In order to get an improved lower bound of eigenfunctions in forbidden regions, we first consider a Neumann problem. 

Let $(\Omega_\Gamma, g)$ be a real-analytic $n$-dimentional Riemannian manifold with smooth boundary $\Gamma$. Now consider following Neumann problem
\begin{align}\label{pde}
(-h^2\Delta_g+ V(x)-E(h) )u_h&=0 \quad \text{in}\, \Omega_\Gamma, \\
\partial_\nu u_h&=0 \quad \text{on}\, \Gamma, \nonumber
\end{align}
here $\partial_\nu$ is the exterior normal derivative, $V(x) >E$ in $\Omega_\Gamma$ and $\| u_h\|_{L^2(\Omega_\Gamma)}\leq 1$.

In the following we fix a small constant $r_0\in(0, \text{inj}(\Omega_\Gamma, g_E) )$ and let $U_\Gamma(r_0)$ be a collar neighbourhood of $\Gamma$ where we have Fermi coordinates $(x', x_n)$ with respect to the Agmon metric $g_E$. The defining function $x_n : \Omega_\Gamma \to \mathbb{R}$ is the distance to the boundary, with the property that $0\leq x_n\leq r_0 $ in $\Omega_\Gamma$ and $x_n = 0$ on $\Gamma$. And $x'$ is constant on geodesics normal to the boundary. 

Denote level set \[\Gamma_\rho= \{x=(x',x_n);\, x_n=\rho, 0\leq \rho\leq r_0\}\]
and let $\gamma_\rho: C^\infty(\Omega_\Gamma) \to C^\infty(\Gamma_\rho)$ be the restriction operator. 

Motivated by \cite{SU} and \cite{GT}, we use the Poisson representation for problem \eqref{pde}, whose parametirx is an Fourier integral operator with complex phase. As the upper half plane model (Sect. \ref{sec2.1}) reveals, the main difficulty of improving the decay rate is to show that the restricted eigenfunctions $u_h\, \vline_\Gamma$ do concentrate near the zero section, indeed we can prove that it localized to frequencies $\geq \lambda^{1/2}h^{1/2}$ is $O(\lambda^{-1/2})$ (Proposition \ref{exterior}). Hence, with establishing a proper lower bound estimate (Proposition \ref{elliptic}), we can get following main theorem of this paper (see Sect. \ref{sec2.2} for the proof).


 

\begin{maintheo}\label{main}
If $u_h$ solves \eqref{pde}, for
$\alpha\in \mathbb{N}^n$, there exist constants $h_0(\alpha)>0$, $\rho_0>0$ and $C(\alpha)>0$ such that for $h\in(0, h_0(\alpha)]$
\begin{equation}\label{result}
\|\partial^\alpha u_h\|_{L^2(\Gamma_\rho)}\geq C(\alpha) e^{-\rho /h} \| u_h\|_{L^2( \Gamma )} \quad \text{if} \,\,  0<\rho< \rho_0.
\end{equation}
\end{maintheo}\

\begin{rem}
Notice that on compact set $\Gamma_\rho$, the $L^2$ norm in ambient metric $g$ and the one in Agmon metric $g_E$ are comparable.
In the following, we shall use $\| \cdot \|_{L^2(\cdot,\, g)}$ to emphasize the $L^2$ norm in ambient metric $g$ instead $\| \cdot \|_{L^2(\cdot)}$ in the conformal Agmon metric $g_E$. Notice 
\[\sqrt{|\det g_E|} = \big( V(x)-E \big)_+^{n/2}\sqrt{|\det g|}.\]
We also use $x=(x',x_n)$ to denote Fermi coordinates in a neighborhood of the boundary $\Gamma$ in ambient metric $g$. For any $U \subset\subset \Gamma_\rho$ in one coordinate patch, by definition
\[ \|\cdot \|^2_{L^2(U,\, g)} = \int_U | \cdot |^2 \, \sqrt{ |\det g_{(x', \rho)} |} \, dx' = \int_U | \cdot |^2 \, \big( V(x', \rho)-E \big)_+^{\frac{n}2} \sqrt{ |\det g_{E(x',\rho)} |} \, d x'.\]
So by partition of unity, for any $W \subset\subset \Gamma_\rho$, $\|\cdot \|_{L^2(W)}$ and $\|\cdot \|_{L^2(W,\, g_E)}$ are comparable, which is equivalent to say, there exist $c, \, C>0$ which are independent of $h$ such that
\begin{equation}
c\|\cdot \|_{L^2(W,\, g)} \leq \|\cdot \|_{L^2(W)} \leq C \|\cdot \|_{L^2(W, \,g)}.
\end{equation}

Hence we can write \eqref{result} in ambient metric
\begin{equation}
\|\partial^\alpha u_h\|_{L^2(\Gamma_\rho,\, g)}\geq C(\alpha) e^{-\rho /h} \| u_h\|_{L^2( \Gamma, \, g )} \quad \text{if} \,\,  0<\rho< \rho_0,
\end{equation}
$\rho$ is the distance from hypersurface $\Gamma_{\rho}$ to hypersurface $\Gamma$ in Agmon metric $g_E$.
\end{rem}

\begin{rem} Back to our problem \eqref{Schrodinger}, if there exists a smooth separating hypersurface $\Gamma$ which is isotopic in classically forbidden region $\Omega_E$  to boundary caustic $\Lambda_E$ and $\Gamma$ bounds a domain $\Omega_\Gamma$ satisfying $\Omega_\Gamma \subset \Omega_E$ (see Figure 1) such that a sequence eigenfunctions $\{u_h\}$ satisfying $\partial_{x_n} u_h \vline_{\Gamma }=0$, then Theorem \ref{main} exactly shows that the sequence eigenfunctions $\{u_h\}$ exponentially decay in Agmon distance from the hypersurface $\Gamma$.
\end{rem}

\begin{figure}[ht]     
\begin{center}
\begin{tikzpicture}[scale=0.9]
\draw [thick, fill=black!10] (-0.2,0) to [out=87, in=170] (2.3, 2.2) to [out=-10, in=5] (4.05, -2.2) to [out=185, in=-15] (2, -1.2) to [out=165, in=-85] (-0.2,0);

\draw [thick] (-1,0) to [out=87, in=170] (2.3, 2.8) to [out=-10, in=5] (4.3, -2.8) to [out=185, in=-15] (1.8, -1.7) to [out=165, in=-85] (-1,0);


\draw [thick, dotted] (0,0) to [out=87, in=180] (2, 2) to [out=0, in=0] (4, -2) to [out=180, in=-15] (2, -1) to [out=165, in=-87] (0,0);

\draw [thick, -latex] (0.8, 0.2) -- (0.23,1);
\node at (-0.3, -0.5) {$\Gamma$};
\node at (-1.3, 0){$\Gamma_\rho$};
\node at (1, 0.2) {$\Lambda_E$};
\node at (5.5,1.5) {$\Omega_\Gamma\subset \Omega_E$};

\end{tikzpicture}
\caption{}
\end{center}
\end{figure}\

At present, we are unable to prove that \eqref{result} holds in the general setting without assuming Neumann condition, but we hope to return to this point elsewhere.

\subsection{Outline of the paper}
In Sect. \ref{sec2.1} we discuss the exponential $L^2$ lower bound for the eigenfunctions $\{u_h\}$ in the upper half plane model. And in Sect. \ref{sec2.2},  we give the proof of Theorem \ref{main}. Finally the key ingredient of the proof, exterior mass estimate (Prop. \ref{exterior}), will be discussed in Sect. \ref{sec3}.\\


\noindent{\sc Acknowledgements:} I would like to thank John Toth for bringing this problem to me and many helpful suggestions at beginning of preparing this manuscript. Many thanks also to Long Jin for all valuable discussions during working on this problem.

Most of this work was done when I was visiting Yau Mathematical Sciences Center, Tsinghua. I am grateful to their kind hospitality and generous financial support.

\bigskip

\bigskip
\section{Exponential lower bound of Analytic Schr\"odinger Eigenfunctions}
In local coordinates, the {\em Laplace-Beltrami} operator $\Delta_g$ has following form,
\[\Delta_g=\sum_{i,j}\frac{1}{\sqrt{|\det g|}} \frac{\partial}{\partial x_j}\left(\sqrt{|\det g|} g^{ij} \frac{\partial}{\partial x_i}\right) ,\]
where $\det g$ is the determinant of $g$.

In local coordinates where $\{V(x)>E\}$, for sufficiently small $h$, the conformal Laplacian under Agmon metric $g_{E}$ is of the form
\[\Delta_{g_{E}}=\frac{1}{V-E} \Delta_g+\sum_{i,j}\frac{n-2}2\frac{ g_E^{ij}}{V-E }{\partial_{x_j} V} \,\partial_{x_i}, \]
here $n\geq 2$ is the dimension of $M$, correspondingly
 \[\frac{P(h)}{V-E}=-h^2\Delta_{g_{E}}+h^2 \frac{n-2}2 \sum_{i,j}\frac{g_E^{ij}}{V-E(h)}{\partial_{x_j}V} \,\partial_{x_i}+1+f(x;h),\]
 here $f(x;h)=\frac{E-E(h)}{V(x)-E(h)}\to 0$ as $h\to 0^+$

In terms of the Fermi coordinates $x=(x', x_n)$, 
\[\frac{P(h)}{V-E}=h^2D_{x_n}^2+ R_h(x, D_{x'})+1+f(x;h) ,\]
where $R_h(x, D_{x'})$ is a second order elliptic differential operator in $x'$ with positive principal symbol $r(x, \xi')$. Indeed $R_h(x, D_{x'})=-h^2\Delta_{g_{E}}\vline_{\,\Gamma_{x_n}}+ h^2\frac{\partial_{x_n} |\det g_{E}| }{2|\det g_{E}|} \partial_{x_n} +h^2\, \frac{n-2}2 \sum_{i,j}\frac{g_{E}^{ij}}{V-E}{\partial_{x_j}V}\partial_{x_i}$, here $\Delta_{g_{E}}\vline_{\,\Gamma_{x_n}}$ is the induced tangential Laplacian, $\Delta_{g_{E}}$, on $\Gamma_{x_n}$.

Set a conjugated operator
\[Q(h)=e^{x_n/h} \frac{P(h)}{V-E} e^{-x_n/h}=h^2D_{x_n}^2+ \tilde R_h(x, D_{x'}) + 2i h D_{x_n}+f(x;h), \]
where $\tilde R_h(x, D_{x'})=e^{x_n/h} R_h(x, D_{x'}) e^{-x_n/h}$ is a second order elliptic differential operator in $x'$ with same positive principal symbol $r(x, \xi')$.



Setting $v_h=e^{\frac{x_n} h} u_h$, then
\begin{align}  \label{rpde}
Q(h)v_h&=0 \quad \,\,\text{in} \,\,\,  \Omega_{\Gamma },
\end{align}
and $v_h(x',0)=u_h(x',0):=\varphi_h(x')$.

Let $U\subset T^*M$ be open. Following \cite{Sj96}, we define the notion of a {\em classical analytic symbol (cl.a.s) of order k} and write $a\in S_{cla}^{m, k}(U)$ provided $a\sim h^{-m}(a_0+ha_1+ \dots)$ in the sense that
\[\partial_x^k \partial_\xi^l \bar\partial_{(x,\xi)} a=O_{k,l}(1)e^{-\left<\xi\right>/Ch}, \quad (x,\xi)\in U,\]
\[\left| a-h^{-m}\sum_{0\leq j\leq \left<\xi\right>/C_0h}h^j a_j\right|=O(1)e^{-\left<\xi\right>/C_1h}, \quad |a_j|\leq C_0C^j j! \left<\xi\right>^{k-j}, \quad (x,\xi)\in U.\]

To keep track of powers of $h$ in the remainders, we use a special class of symbols than one used in \cite{DJ}. 

Fix parameter $0\leq \rho <1$, we say that an $h$-dependent symbol $a$ lies in the class $S^{ \text{comp}}_{\rho}(U)$ if

(1) $a(x,\xi ;h)$ is smooth in $(x,\xi)$ in $U$, defined for $0<h\leq 1$, and supported in an $h$-independent compact subset of $U$;

(2) $a$ satisfies the derivative bounds
\begin{equation}
\sup_{x,\xi}  |\partial_x^{\alpha} \partial_\xi^{\beta} a(x,\xi; h)|\leq C_{\alpha \beta} h^{ -\rho|\beta|}.
\end{equation}

\begin{rem}
$S^{ \text{comp}}_{\rho}(U)$ is a special class of $S^{ \text{comp}}_{L, \rho,\rho'}(U)$ used in \cite[Appendix]{DJ}. As a model case, one can take $U=T^*\mathbb{R}^n$, $L=L_0=\text{span}\{\partial_{x_1},\dots, \partial_{x_n}\}$ and $\rho'=0$.
\end{rem}

The corresponding standard semiclassical pseudodifferential operators have Schwartz kernels that are sums of the local integrals of the form
\[Op_h(a) (x,y) :=\frac{1}{ (2\pi h)^n } \int_{\mathbb{R}^n } e^{\frac{i}h \left<x-y, \xi \right> } a(x,\xi ; h)  d\xi . \]
\medskip

\subsection{Upper half plane model}\label{sec2.1}
Let us first explain heuristically where the proper h-exponential decay comes from. Consider $L^2$-normalized solutions $\{u_h\}$ on the upper half plane $\mathbb{R}^n_+=\{(x',x_n)\in \mathbb{R}^n;\, x_n>0\}$,
\begin{align}\label{half}
\left(-h^2\Delta+1 \right)u_h(x)&=0 \quad \text{on}\, \mathbb{R}^n_+,  \nonumber \\
u_h(x',0)&=\phi_h(x') \quad \text{on}\, \partial \mathbb{R}^n_+.
\end{align}
We assume that $\phi_h(x')$ concentrates near zero section $\xi'=0$. More explicitly, 
\begin{equation}\label{exterior ex}
\|Op_h(1-\chi_\delta(\xi'))\varphi_h\|_{L^2(\partial\mathbb{R}^n_+)}\leq \epsilon \|\varphi_h\|_{L^2(\partial\mathbb{R}^n_+)},
\end{equation}
here $\epsilon$ is a positively small constant and $\chi_\delta\in C_0^\infty(T^*\mathbb{R}^n_+; [0,1])$ is a cutoff supported near the zero section, with $\chi_\delta(x',\xi')=1$ for $\{0\leq |\xi'|\leq \delta/2\}$ and $\chi_\delta(x',\xi')=0$ for $|\xi'|>\delta$. Here $\delta>0$ is some arbitrarily small but fixed constant.

It's straightforward to check that the restricted Poisson operator of \eqref{half} is of the form 
\[\gamma_\rho K w(x', \rho)=\frac{1}{(2\pi h)^{n-1}}\int\int e^{\frac{i}h\left<x'-y', \xi'\right>-\frac{\rho}h\sqrt{|\xi'|^2+1}}  w(y')d y' d \xi',\]


Define the {\em semiclassical Fourier transform} for $h>0$
\[\mathcal{F}_hu(\xi):=\int_{\mathbb{R}^n}e^{-\frac{i}h\left<x,\xi\right>} u(x)dx\]
and its inverse
\[\mathcal{F}^{-1}_h v(x):=\frac{1}{(2\pi h)^n}\int_{\mathbb{R}^n}e^{\frac{i}h\left<x,\xi\right>} v(\xi)d\xi .\]
Denote $\partial\mathbb{R}^n_\rho =\{(x', \rho )\in \mathbb{R}^n;\, \rho \text{ is a fixed positive constant}\}$ and apply Plancherel formula
\begin{align}\label{upper half}
\|(\gamma_\rho K) \varphi_h\|_{L^2(\partial \mathbb{R}^n_\rho)}&=\|\mathcal{F}^{-1}_h\big(e^{-\frac{\rho}h\sqrt{|\xi'|^2+1}}(\mathcal{F}_h\varphi_h)\big)\|_{L^2(\mathbb{R}^{n-1})} \nonumber\\
&=\frac{1}{(2\pi h)^{(n-1)/2} }\|e^{-\frac{\rho}h\sqrt{|\xi'|^2+1}}(\mathcal{F}_h\varphi_h)\|_{L^2(\mathbb{R}^{n-1})}  \nonumber \\
&\geq \frac{1}{(2\pi h)^{(n-1)/2} }\|e^{-\frac{\rho}h\sqrt{|\xi'|^2+1}}\chi_\delta(\xi')(\mathcal{F}_h\varphi_h)\|_{L^2(\mathbb{R}^{n-1})} \nonumber \\
&\geq e^{-\frac{\rho}h\sqrt{\delta^2+1}}\frac{1}{(2\pi h)^{(n-1)/2} }\|\chi_\delta(\xi')(\mathcal{F}_h\varphi_h)\|_{L^2(\mathbb{R}^{n-1})} \nonumber \\
&=e^{-\frac{\rho}h\sqrt{\delta^2+1}} \|Op_h(\chi_\delta)\varphi_h\|_{L^2(\partial \mathbb{R}^n_+)}.
\end{align}

With the help of \eqref{exterior ex}, for small $h$, the RHS of \eqref{upper half} is
\begin{align*}
&\geq e^{-\frac{ \rho}h\sqrt{\delta^2+1}}\left ( \|\phi_h\|^2_{L^2(\partial \mathbb{R}^n_+)} - \|Op_h(1-\chi_\delta(\xi'))\varphi_h\|_{L^2(\partial\mathbb{R}^n_+)} \right ) \\
&\geq \frac{1}2 e^{-\frac{ \rho}h\sqrt{\delta^2+1}} \|\phi_h\|^2_{L^2(\partial \mathbb{R}^n_+)} .
\end{align*}
 
Consequently, under assumption \eqref{exterior ex}
\begin{align}\label{upper half result}
\|u_h\|_{L^2(\partial \mathbb{R}^n_{\rho})}&=\left<\gamma_\rho K\varphi_h ,\gamma_\rho K\varphi_h\right>^{1/2}_{\partial \mathbb{R}^n_{\rho}}
\geq \frac{1}2 e^{-\frac{ \rho}h\sqrt{\delta^2+1}} \|\phi_h\|^2_{L^2(\partial \mathbb{R}^n_+)}
\end{align}
for sufficiently small $h>0$.

\begin{rem}
If the boundary data $\varphi_h(x')$ is independent of $h$, we have similar inequality \eqref{upper half result} with directly applying analytic stationary phase theorem without assumption \eqref{exterior ex}.
\end{rem}
From above argument, inequality \eqref{exterior ex} is a key ingredient. Intuitively, in order to get an improved reverse Agmon estimate, we need to establish similar inequality(Proposition \ref{exterior}).
\medskip

\medskip
\subsection{Lower bound of the Poisson parametrix: proof of Theorem \ref{main}} \label{sec2.2}


%
%


Throughout this section, we always work in the collar neighbourhood $U_{\Gamma}(r_0)$. First choose a covering of $\Gamma_\rho$ by finitely many coordinate charts $\{W_j\times \{\rho \} \}_{j\in J}$ and a corresponding smooth, locally finite partition of unity $\{\chi_j\}_{j\in J}$, with
\[\sum_{j\in J} \chi_j(x')=1,\,\, (x', \rho)\in W_j\times \{\rho\} \subset \Gamma_\rho,\,\,  \chi_j(x')\in C^{\infty}_0(W_j). \]
For convenience set $W= \cup_j W_j$, so $\Gamma_\rho=W\times \{\rho\}$.

It follows from \cite{GT} (or \cite{SU}) the Poisson kernel $K(x', x_n, y', \xi')$ of  problem \eqref{rpde} has following form
\begin{align} \label{kernel}
K(x', x_n, y', h)=\frac{1}{(2\pi h)^{n-1}}\sum_{j,k} \chi_j(x')  \int_{\mathbb{R}^{n-1}} e^{\frac{i}h[\phi(x', x_n,\xi')-y'\cdot \xi']}a(x', x_n,\xi'; h) \chi_k(y') d\xi'  \nonumber \\
+O (e^{-C_0/h}), 
\end{align}

\noindent here $a(x', x_n, \xi'; h)$ on $U_\Gamma(r_0)\times \mathbb{R}^{n-1}$ is a cl.a.s of order $0$ and the summation is for all $\chi_j$ and $\chi_k$ such that $\text{supp}(\chi_j) \cap \text{supp}(\chi_k) \neq \emptyset$. Also $a(x', 0,\xi'; h)=1+O(h)$, and there exist constants $\rho_0, C_1, C_2>0$, such that 
\begin{equation}\label{a}
C_1\leq a(x', x_n, \xi'; h)\leq C_2 \quad \text{for}\,\, \forall\,\, 0 \leq x_n\leq \rho_0.
\end{equation}

\noindent Here $\phi$ solves following Hamilton-Jacobi equation
\begin{equation}\label{HJ}
(\partial_{x_n} \phi)^2 + 2i (\partial_{x_n} \phi)+r(x', x_n, \partial_{x'}\varphi ) =0, \quad \phi(x', 0, \xi')=\left<x', \xi'\right>,
\end{equation}
with $r(x', x_n, \xi'):= \Sigma_{i,j}^{n-1} g_{E}^{ij}(x', x_n)\xi_i \xi_j $, the semiclassical principle symbol of $\tilde R_h(x', x_n, D_{x'})$. For simplicity, we write $r(x', x_n, \xi')=  |\xi'|^2_{ (x', x_n) } $. 
Such $\phi$ in \eqref{HJ} exists by applying Cauchy-Kowalevski theorem since polynomial $r(x', x_n, \xi')$ is analytic. 

By the positivity of the metric $g_E$ for $x_n\geq 0$,
\begin{equation} \label{metric positivity}
 \exists\, C, \tilde C>0 \quad \text{such that} \,\, C |\xi'|^2 \leq |\xi'|^2_{ (x', x_n) } \leq \tilde C |\xi'|^2 \, , \quad \text{for} \,\, \forall (x' ,x_n) \in \,\, U_{\Gamma}(r_0),\,\, x_n\geq 0. 
 \end{equation}

With the natural branch of $r^{1/2}$ with a cut along the real negative axis, it follows by Taylor expansion in $\xi'$ near $0$ that


\begin{equation} \label{phase}
\phi= \left<x', \xi'\right>+ i\phi_1(x',x_n,\xi')
\end{equation}
where 
\begin{align} \label{phi_1}
\phi_1(x', x_n, \xi')=\sum_{0\leq |\alpha| \leq k-1} \frac{\partial^\alpha \phi_1 (x', x_n, 0)}{\alpha !} (\xi')^\alpha + R_k, \quad \varphi_1(x', 0, \xi')=0.
\end{align}

For $\alpha=0$, by \eqref{HJ} $\varphi_1(x', x_n, 0)$ solves following equation
\[ (\partial_{x_n}\varphi_1)^2+2\partial_{x_n}\varphi_1- \Sigma_{i,j}^{n-1} g_{E}^{ij}(x', x_n)\partial_{x_i}\varphi_1\partial_{x_j}\varphi_1=0, \quad \varphi_1(x', 0, 0)=0.\]
By the uniqueness of Cauchy-Kowalevski theorem, $\varphi_1(x', x_n, 0)=0$ is the unique solution.
Similarly, we have $\partial_{\xi_j} \varphi_1(x', x_n, 0)=0$ for $1\leq j \leq n-1$. Hence,
\begin{equation}\label{phi1}
\phi_1(x', x_n, \xi')=O_{x', x_n}(|\xi'|^2) \quad  \text{for small}\, |\xi'|.
\end{equation}

\noindent By Taylor expansion near $x_n=0$, we also have
\begin{equation}\label{phi1 xn}
\varphi_1(x', x_n, \xi')= x_n(\sqrt{|\xi'|_{(x',0)}+1}-1)+O\left(x_n^2 |\xi'|_{(x',0)}\right) \quad  \text{for large}\, |\xi'|.
\end{equation}



Modulo $O(e^{-C_0/h})$ term, we write
\[\gamma_\rho K=Op_h(\sigma_{\gamma_\rho K}), \] 
where $\sigma_{\gamma_\rho K}(x', \xi'):=\sum \chi_i e^{-\phi_1(x',\, \rho, \, \xi') / h} a(x',\rho,\xi'; h)$.

Recall that $\varphi_h(x'):=v_h\vline_\Gamma=u_h\vline_\Gamma$. One has
\begin{align}\label{identity}
\gamma_\rho K \varphi_h :&= e^{ \rho /h}u_h(x', \rho)  \nonumber \\
&= \frac{1}{(2\pi h)^{n-1}}\sum_{j,k} \chi_j(x') \int_{W_j} \int_{\mathbb{R}^{n-1}}  e^{\frac{i}h[\phi(x', \rho,\xi')-y'\cdot \xi']}a(x', x_n,\xi'; h) \chi_k(y') \varphi_h(y') d\xi' dy' \nonumber \\
&\qquad\qquad\qquad\qquad\qquad\qquad\qquad\qquad\qquad\qquad +O(e^{-C_0/h})  \varphi_h(x') .
\end{align}

\begin{figure}[ht]
\begin{tikzpicture}
\draw[->] (-0.5,0)--(0,0) node[below left]{\tiny 0}--(1.5, 0) node[below]{\tiny $ M $}  -- (6,0) node[below]{$ x $};
\draw[fill] (0, 1) circle [radius=1pt];
\draw[fill] (1.5,0) circle [radius=1pt];
\draw[dashed] (0, -0.3) -- (0, 1.3);
\draw[thick] (0, 1) to [out=-45, in=180] (1.5, 0.3) to [out=0, in=180] (5.5, 0.3);
\node[above] at (3, 0.4) {$b( \cdot )$};
\node[right] at (0, 1.2) {\tiny 1};
\node[right] at (5.5, 0.5) {\tiny $\zeta$};
\end{tikzpicture}
\caption{}
\end{figure} 

Choose a decreasing smooth function $b(x)\geq e^{-x}$ for $x\geq 0$  (see Figure 2) satisfying that
\begin{align*}
b(x) &= e^{-x}, \quad 0\leq x \leq \frac{M} {2}  , \\
b(x) &=\zeta, \quad  x \geq M, 
\end{align*}
and $|\frac{d^n}{dx^n} b(x)|\leq C_n b(x)$, here $M$ is a large constant which will be determined later and $e^{-M}< \zeta<e^{-\frac{M}{2} }$ is a constant.

Now set
\begin{align*} b_j( x', \rho, \xi' ; h)&= \chi_j(x')\, b(\varphi_1(x', \rho, \xi')/h ).
\end{align*}
By the definition of $b$ and the positivity \eqref{metric positivity} of the metric, one has 
\begin{equation} \label{xi}
b_j( x', \rho, \xi' ; h)\equiv \chi_j e^{-\varphi_1(x', \rho, \xi')/h} \quad \text{if}\,\,  |\xi'|^2 \leq cM h, 
\end{equation}
and
\begin{equation}
b_j( x', \rho, \xi' ; h) \equiv \chi_j  \zeta \quad \text{if}\,\,  |\xi'|^2 \geq c' M h, 
\end{equation}
here $c'>c$ are two positive constants. 

One also has
\begin{equation}
b_j(x' , \rho, \xi'; h) \geq  \chi_j \zeta, \quad \text{for}\,\,  0 \leq \rho < \rho_0.
\end{equation}
Using \eqref{phi1}, this is straightforward to check that for $|\xi'|^2 \leq c' M h$
\begin{equation*}
|\partial^\alpha_{x'} \varphi_1 | \leq C_\alpha  M h \,\, \text{if} \,\, |\alpha|\geq 1,
\end{equation*}
\begin{equation} \label{d varphi}
|\partial^\beta_{\xi'} \varphi_1 | \leq C_\beta M^{1/2}  h^{ 1/2 } \,\, \text{if,} \,\, |\beta|=1 \quad  |\partial^\beta_{\xi'}  \varphi_1 | \leq C_\beta \,\, \text{if} \,\, |\beta| \geq 2.
\end{equation}

By chain rule, this is straightforward to check that
\begin{align}
| \partial_{x_k} b_j | \leq  C_M  b_j, \quad | \partial_{\xi_k} b_j | \leq  C_M h^{-1/2} b_j, \nonumber \\
 | \partial_{x_k x_l} b_j | \leq  C_M b_j , \quad | \partial_{\xi_k \xi_l} b_j | \leq  C_M h^{-1} b_j .
\end{align} 
Furthermore by {\em Leibniz} rule, one has
\begin{equation}\label{b_j}
|\partial_{x'}^\alpha \partial_{\xi'}^\beta b_j | \leq C_{\alpha \beta}  h^{ - \frac{|\beta|}2 } b_j,\end{equation} 
which means that $ b_j(x' , \rho, \xi'; h)-\chi_j \zeta \in S^{\text{comp}}_{\frac{1}2 }(T^* W_j) $.

For convenience, set $f_\rho(x', \xi'; h):= \sum_j  b_j(x' , \rho, \xi'; h)$. From \eqref{b_j}, one has
\begin{equation}\label{frho}
|\partial_{x'}^\alpha \partial_{\xi'}^\beta f_\rho | \leq C_{\alpha \beta}  h^{-\frac{|\beta|}2 } f_\rho.
\end{equation}
which means that $f_\rho(x', \xi'; h)- \zeta\in S^{\text{comp}}_{\frac{1}2 } (T^* W) $.

\begin{prop} \label{elliptic}
There exist $h_0>0$ and $C_3>0$ such that
\begin{equation}\label{positivity}
 \|  Op_h(f_\rho \, a)  \phi_h \|_{L^2(\Gamma_\rho)} \geq  C_3 \zeta \| \phi_h\|_{L^2(W) }.
\end{equation}
for $0<h< h_0$.
\end{prop}
\begin{proof}



The proof essentially follows \cite[Theorem 4.29]{Zwo}. From the definition of $f_\rho$ and \eqref{frho}, we know $f_\rho^{-1} -\zeta^{-1} \in f_\rho^{-1} S^{\text{comp}}_{\frac{1}2 } (T^* W)$.

By composition formula (see \cite[(4.3.16)]{Zwo} and \cite[(A.23)]{DJ}) and \eqref{frho}, one can get
\begin{equation}
f_\rho \# {f_\rho }^{-1} (x', \xi' ; h)=1  +r_1
\end{equation}
with
\begin{equation*}
r_1\in  h^{\frac{1}2} S^{\text{comp}}_{\frac{1}2 } (T^* W)
\end{equation*}
Likewise
\begin{equation}
f_\rho^{-1}\# f_\rho (x', \xi'; h) =1+r_2 \quad \text{with} \,\,  r_2 \in h^{\frac{1}2} S^{\text{comp}}_{\frac{1}2 } (T^* W).
\end{equation}
Hence if $R_1:= Op_h(r_1)$ and $R_2:= Op_h(r_2)$, we have
\[Op_h(f_\rho)\, Op_h(f_\rho^{-1}) = I+R_1,\]
\[Op_h(f_\rho^{-1})\, Op_h(f_\rho) = I+R_2,\]
with $\|R_1\|_{L^2\to L^2}, \|R_2\|_{L^2\to L^2}= O(h^{\frac{1}2})  \leq 1/2$.

It remains to apply triangle inequality and \cite[Lemma A.5]{DJ} to get that
\begin{equation}
 \| Op_h( f_\rho^{-1}) \|_{L^2\to L^2} \leq  \| Op_h( f_\rho^{-1}-\zeta^{-1}) \|_{L^2\to L^2} +\zeta^{-1}\leq C \zeta^{-1},
\end{equation}
with noticing that $1\leq f_\rho^{-1} \leq \zeta^{-1}$ and $f_\rho^{-1} -\zeta^{-1} \in f_\rho^{-1} S^{\text{comp}}_{\frac{1}2 } (T^* W)$.

Then with the help of \eqref{a},
\begin{align*}
\frac{C_1}2 \| \phi_h\|_{L^2(W) } &\leq  \| Op_h(a) \phi_h\|_{L^2(\Gamma_\rho) } \\
&\leq \|(I+R_2)^{-1} Op_h(f_\rho^{-1}) Op_h(f_\rho) Op_h(a) \phi_h\|_{L^2(\Gamma_\rho) } \\
&\leq  C \zeta^{-1} \| Op_h(f_\rho) Op_h(a) \phi_h\|_{L^2(\Gamma_\rho) } \\
&\leq  \frac{C}{2} \zeta^{-1} \| Op_h(f_\rho a) \phi_h\|_{L^2(\Gamma_\rho) }
\end{align*}
implies the result. The last line follows by applying composition formula.

\end{proof}


\begin{figure}[ht]
\begin{tikzpicture}
\draw[->] (-0.5,0)--(0,0) node[below left]{\tiny 0 } --(1, 0) node[below]{\tiny $ O(M^{\frac{1}2 } h^{\frac{1}2 })$}   -- (5.5,0) node[below]{\tiny $ |\xi'|  $};
\draw[dashed] (-0.5, 1) -- (5, 1);
\draw[fill] (0,1) circle [radius=1pt];
\draw[fill] (1,0) circle [radius=0.7 pt];
\draw[dashed] (0, -0.3) -- (0, 1.7);
\draw[thick] (1, 0) to [out=0, in=180] (5, 0.9);
\node[right] at (0,1.25) {\tiny $\zeta a$};
\node[right] at (5, 0.75) {$f_\rho a  -\sigma_{\gamma_{\rho K} }$};
\end{tikzpicture}
\caption{}
\end{figure}


\begin{prop} \label{bound}
$Op_h( f_\rho a -\sigma_{\gamma_\rho K}): L^2( \Gamma ) \to L^2( \Gamma_\rho )$ is a bounded operator. More explicitly, there exists a constant $ C_4>0$ such that 
\[ \|Op_h( f_\rho a -\sigma_{\gamma_\rho K})\|_{L^2( \Gamma  ) \to L^2( \Gamma_\rho )} \leq C_4 \zeta. \]
\end{prop}
\begin{proof}


From the definition of $f_\rho$, one has $ f_\rho - \sum \chi_i e^{-\phi_1(x',\, \rho, \, \xi') / h}   \leq \zeta$ which is a {\em Schwartz} function in $\xi'$ (Figure 3).

Then by $L^2$ boundedness, one can conclude that 
\[ \|Op_h( f_\rho a -\sigma_{\gamma_\rho K})\|_{L^2( \Gamma  ) \to L^2( \Gamma_\rho )} = \|Op_h( f_\rho a - \sum \chi_i e^{-\phi_1(x',\, \rho, \, \xi') / h} a) \|_{L^2( \Gamma  ) \to L^2( \Gamma_\rho )} \leq C_4 \zeta. \]

\end{proof}

Taking radial cutoffs $\chi_{\text{in}}, \, \chi_{\text{out}} \in C^\infty(\mathbb{R}; [0,1])$ with $\chi_{\text{in}}(\xi')=\chi_{\text{in}}(|\xi'|; h, M ) $ and $\chi_{\text{out}}(\xi')=\chi_{\text{out}}(|\xi'|; h, M ) $ such that \[\chi_{\text{in}}(\xi') + \chi_{\text{out}}(\xi')=1 \] and
\begin{align*}
\chi_{\text{in}} (\xi' ) &= 1, \quad |\xi'|^2 \leq \frac{c}2 M h , \\
\chi_{\text{in}} (\xi') &= 0, \quad  |\xi'|^2 \geq c M h, 
\end{align*}
here $c$ is a constant in \eqref{xi}.

According to \eqref{xi}, one has
\[ (f_\rho a -\sigma_{\gamma_\rho K} ) \chi_{\text{out}} =f_\rho a - \sigma_{\gamma_\rho K}.\]

\noindent By Proposition \ref{bound} 
\begin{align} \label{ineq 2}
& \| Op_h(f_\rho \,a-\sigma_{\gamma_\rho K}) \phi_h\|_{L^2( \Gamma_\rho )} \nonumber\\
= & \| Op_h\Big( (f_\rho \,a-\sigma_{\gamma_\rho K})  \chi_{\text{out}} \Big)  \phi_h\|_{L^2( \Gamma_\rho )}  \nonumber\\
= &  \|   Op_h(f_\rho \,a-\sigma_{\gamma_\rho K}) Op_h (\chi_{\text{out}} ) \phi_h\|_{L^2( \Gamma_\rho )}  \nonumber \\
\leq & C_4 \zeta \|Op_h ( \chi_{\text{out}}) \phi_h\|_{L^2(W )},
\end{align}
in the 2nd last line, we applied the composition formula of operators with noticing that $\chi_{\text{out}} $ only depends on $\xi'$.


Thereafter, by triangle inequality and along with \eqref{positivity} and \eqref{ineq 2}, from \eqref{identity}  one can conclude that
\begin{align} \label{main eq1}
 \| e^{ \rho  /h}u_h(x', \rho) \|_{L^2( \Gamma_\rho )}
= & \|Op_h \big(f_\rho \,a- (f_\rho \,a-\sigma_{\gamma_\rho K}) \big) \varphi_h \|_{L^2( \Gamma_\rho )} -  O(e^{-C_0/h})\| \phi_h\|_{L^2(W  )} \nonumber \\
\geq &\|Op_h(f_\rho \,a ) \varphi_h \|_{L^2( \Gamma_\rho )}-  \|Op_h(f_\rho \,a-\sigma_{\gamma_\rho K}) \varphi_h \|_{L^2( \Gamma_\rho )}  \nonumber \\
& -  O(e^{-C_0/h})\| \phi_h\|_{L^2(W  )} \nonumber \\
\geq & C_3 \zeta  \| \phi_h\|_{L^2(W  )}- C_4 \zeta \|Op_h(\chi_{\text{out}} ) \phi_h\|_{L^2( W  )} \nonumber \\
& -  O(e^{-C_0/h})\| \phi_h\|_{L^2(W  )}.
\end{align}

From Proposition \ref{exterior}, with taking $h$ sufficiently small and $M$ sufficiently large, one has
\begin{equation}\label{main eq2}
C_4 \zeta \|Op_h(\chi_{\text{out}}(|\xi'|; h, M) ) \phi_h\|_{L^2( W  )} << C_3 \zeta  \| \phi_h\|_{L^2( W )}.
\end{equation}

Consequently, from \eqref{main eq1} and \eqref{main eq2}, one has
\begin{equation} \label{lower1}
\| u_h \|_{L^2( \Gamma_\rho )} \geq  C e^{-\rho /h} \| \phi_h\|_{L^2( W  )} =C e^{-\rho /h} \| u_h\|_{L^2( \Gamma  )} \qquad \text{if}\,\, \rho \leq \rho_0. 
\end{equation}




The same argument works for $\| \partial^\alpha u_h\|^2_{L^2(\Gamma_{\rho})}$ with each differentiation creating a power of $h^{-1}$.

\bigskip

\section{Exterior mass estimates on $\Gamma$}\label{sec3}
Before presenting our main result in this section, we introduce following lemma.

\begin{lem}\label{lem1}
For small $\rho> 0$ one has
\begin{equation} \label{upper1}
\|u_h\|_{L^2( \Gamma_ \rho)} \leq C(\rho) e^{- \rho /h} \| u_h(x) \|_{L^2( \Gamma )}.
\end{equation}
\end{lem}
\begin{proof}
From equation \eqref{rpde} and \eqref{kernel}, we know for any $(x', \rho) \in \Gamma_\rho$
\begin{align} \label{u}
e^{ \rho /h} &u_h(x', \rho ) \nonumber \\
&=\frac{1}{(2\pi h)^{n-1}}\int_{W} \int_{\mathbb{R}^{n-1}}  e^{\frac{i}h[\phi(x', \rho ,\xi')-y'\cdot \xi']}a(x',\rho ,\xi'; h) \chi(x'-y') \varphi_h (y') d\xi' dy'  \nonumber \\
&\qquad\qquad\qquad\qquad\qquad\qquad\qquad\qquad\qquad\qquad\qquad\qquad +O(e^{-C_0/h}) \varphi_h(x') \nonumber \\
&=\frac{1}{(2\pi h)^{n-1}}\int_{W} \int_{|\xi'|\leq R}  e^{\frac{i}h[\phi(x', \rho ,\xi')-y'\cdot \xi']}a(x', \rho ,\xi'; h)\chi(x'-y') \varphi_h (y') d\xi' dy' \nonumber \\
&\quad + \frac{1}{(2\pi h)^{n-1}}\int_{W} \int_{|\xi'|\geq R }  e^{\frac{i}h[\phi(x', \rho ,\xi')-y'\cdot \xi']}a(x', \rho ,\xi'; h) \chi(x'-y') \varphi_h (y') d\xi' dy'  \nonumber \\
&\qquad\qquad\qquad\qquad\qquad\qquad\qquad\qquad\qquad\qquad\qquad\qquad +O(e^{-C_0/h}) \varphi_h(x') ,
\end{align}
where $\chi(x'-y')$ is a near-diagonal cutoff.
With help of equation \eqref{phi1 xn}, one can apply Cauchy-Schwartz inequality to get that
the integrations in \eqref{u} are bounded by $C(\rho) \|\varphi_h \|_{L^2}$. Then we can complete the proof.

\end{proof}
\begin{rem}
Like Agmon estimate, this lemma shows how fast the restricted eigenfunctions decay from a hypersurface in Forbidden regions.
From \eqref{upper1}, we know that \eqref{lower1} is sharp.
\end{rem}




From now on we fix a small constant $ r_{0,g}\in(0, \text{inj}(\Omega_\Gamma, g) )$ and let $U_\Gamma(r_{0, g} )$ be a collar neighbourhood of $\Gamma$ where we have Fermi coordinates $(x', x_n)$ with respect to the ambient metric $g$. We also set \[\Gamma_{\rho,g}= \{x=(x',x_n);\, x_n=\rho, 0\leq \rho\leq r_{0,g}\}.\]

\begin{lem}\label{lem2}
For small $\rho>0$ one has
\begin{equation}
\|\partial_{x_n} u_h\|_{L^2(\Gamma_{\rho, g})} \leq C(\rho) h^{-1} \|u_h\|_{L^2(\Gamma)}
\end{equation}
\end{lem}
\begin{proof}
Like we did in section \ref{sec2.2}, in ambient metric $g$ one has
\begin{align}
u_h(x', x_n)=\tilde K \varphi_h+ O(e^{-C_0/h}\varphi_h),
\end{align}
where
\begin{align}
\tilde K(x', x_n, y', h)=\frac{1}{(2\pi h)^{n-1}}\sum_{j,k} \chi_j(x')  \int_{\mathbb{R}^{n-1}} e^{\frac{i}h[\tilde\phi(x', x_n,\xi')-y'\cdot \xi']} \tilde a(x', x_n,\xi'; h) \chi_k(y') d\xi'  \nonumber \\
+O (e^{-C_0/h}), 
\end{align}

\noindent here $\tilde a(x', x_n, \xi'; h)$ on $U_\Gamma(r_{0,g})\times \mathbb{R}^{n-1}$ is a cl.a.s of order $0$ and $\tilde a(x', 0,\xi'; h)=1+O(h)$.

 $\tilde\phi$ solves following Hamilton-Jacobi equation
\begin{equation}
(\partial_{x_n} \tilde \phi)^2 + V -E + \tilde r(x', x_n, \partial_{x'}\tilde\varphi ) =0, \quad \tilde\phi(x', 0, \xi')=\left<x', \xi'\right>,
\end{equation}
with $\tilde r(x', x_n, \xi'):= \Sigma_{i,j}^{n-1} g^{ij}(x', x_n)\xi_i \xi_j $. For simplicity, we write $\tilde r(x', x_n, \xi')=  |\xi'|^2_{ (x', x_n) } $. 

With the natural branch of $\tilde r^{1/2}$ with a cut along the real negative axis, it follows by Taylor expansion in $x_n$ near $0$ that \[\tilde \varphi=\left<x', \xi' \right> +i \tilde \varphi_1(x', x_n, \xi'),\]
where
\begin{equation}
\tilde\varphi_1(x', x_n, \xi')= x_n\left(\sqrt{V(x',0)-E + |\xi'|^2_{(x',0)} }\right)+O\left(x_n^2 |\xi'|_{(x',0)}\right) \quad  \text{for large}\, |\xi'|.
\end{equation}
Then we can finish the proof following the steps in Lemma \ref{lem1} with noticing that each differentiation creating a power of $h^{-1}$.

\end{proof}

Following results essentially follow \cite{BHT}. To emphasize the normal direction $x_n$, let us write
\begin{equation}\label{change variable}
r= x_n,
\end{equation}
\noindent and set $P(r)=-h^2\Delta_{g}\vline_{\,\Gamma_r}$.

Let the smooth function $\chi: \mathbb{R} \to [0,1]$ satisfy $\chi(t)=0$ for $t\leq 1$ and $\chi(t)=1$ for $t\geq 2$. Define $f_\lambda(t)$ by
\[f_\lambda(t)=\chi \left( \frac{t}{ \lambda h } \right), \]
here $\lambda$ is a large parameter.

\begin{lem}\label{lemma}
One has 
\[ \|(d/ dr)^k f_\lambda(P)\|_{L^2( \Gamma_{r,g} )} \leq C_k  \lambda^{-k} h^{-k }\]
for $k=0,\, 1,\, 2,$ and $h\leq 1$.
\end{lem}
For the proof, one can refer \cite[Lemma 3.3 and Lemma 5.1]{BHT}.

\begin{proof}
We select an almost analytic extension of $f_\lambda$, denoted by $F_\lambda$. We extend the smooth function $\chi (t)$ to an almost analytic extension $\tilde \chi(z)$, satisfying
\[ |\partial_z^k \bar\partial \tilde \chi(z)|\leq C_{k, N}|Im\, z|^N \quad \forall z\in \mathbb{C}, \,\forall N\in \mathbb{N},\]
here $\bar\partial$ is the usual $d$-bar operator, $\partial_x+i \partial_y$.
We then define
\[F_\lambda (z)= \tilde \chi \left( \frac{z}{ \lambda h } \right).\]
Note that $F_\lambda (t)= f_\lambda (t)$ for real $t$, and since $P(r)\geq 0$ as an operator, it follows that $F_\lambda(P)= f_\lambda(P)$. Due to the scaling in the definition of $h$, we have
\begin{equation}
|\partial_z^k \bar\partial F_\lambda(z)|\leq C_{k, N} \lambda^{-(k+N+1)} h^{-(k+N+1)} |Im\, z|^N \quad \forall z\in \mathbb{C}, \,\forall N\in \mathbb{N}.
\end{equation}\label{almost analytic}
In addition, we can assume that $\bar\partial \tilde\chi$ is supported in the set $[1, 2] \times i[-1, 1]$. Consequently, $\bar\partial F_\lambda$ is supported in $[\lambda h , 2 \lambda h ]\times i[- \lambda h , \lambda h] $, which is a set of measure $O(\lambda^2 h^{2} )$.

We can express $f_\lambda(P)$ in terms of $F_\lambda$ using the standard Helffer-Sj\"ostrand formula (see Theorem 8.1 \cite{Zwo})
\begin{equation}
f_\lambda(P) = \frac{1}{2\pi}\int_{\mathbb{C}} \bar\partial F_\lambda(z) (P-z)^{-1}d L(z).
\end{equation}
Here the integral is over the entire complex plane $\mathbb{C}$, $L(z)$ is Lebesgue measure on $\mathbb{C}$. Using this formula we can easily express the $r$-derivatives of $f_\lambda(P)$. For example, using
\[\dot{(P-z)}^{-1}=-(P-z)^{-1}\dot P (P-z)^{-1},\]
with dots indicating differentiation with respect to $r$.
We have
\begin{align}
\dot{f_\lambda(P)} =&-\frac{1}{2\pi}\int_{\mathbb{C}} \bar\partial F_\lambda(z) (P-z)^{-1}\dot P (P-z)^{-1} d L(z), \nonumber \\
\ddot{f_\lambda(P)} =& \frac{2}{2\pi} \int_{\mathbb{C}} \bar\partial F_\lambda(z) (P-z)^{-1}\dot P (P-z)^{-1} \dot P (P-z)^{-1} d L(z) \\
&- \frac{1}{2\pi} \int_{\mathbb{C}} \bar\partial F_\lambda(z) (P-z)^{-1}\ddot P (P-z)^{-1} d L(z). \nonumber
\end{align}


Note that $1+P$ is invertible, so we have
\[ (P-z)^{-1}\dot P (P-z)^{-1}=  (P-z)^{-1}(1+P) (1+P)^{-1} \dot P (1+P)^{-1} (1+P) (P-z)^{-1}.\]

Note that $(1+P)^{-1} \dot P (1+P)^{-1} $ has an $O(1)$ operator norm bound on $L^2(\Gamma_{r,g} )$ uniform in $r$ and $h$. On the other hand, using spectral theory, the operator norms of $(1+P) (P-z)^{-1}$ and $(P-z)^{-1} (1+P)$ are, for $z\in \supp \bar \partial F_\lambda $, bounded by
\[\sup_{t\in[0, +\infty) }| (1+t)(t-z)^{-1} | \leq C |\Im z|^{-1}.\]
Hence $(P-z)^{-1}\dot P (P-z)^{-1}$ is bounded by $C |\Im z|^{-2}$ for $z\in \supp \bar \partial F_\lambda $.

Together with the bound on $\bar\partial F_\lambda(z)$
\[ | \bar\partial F_\lambda (z)|\leq C_{N} \lambda^{-(N+1)} h^{-(N+1)} |Im\, z|^N \quad \forall z\in \mathbb{C}, \,\forall N\in \mathbb{N},\]
along with the $O(\lambda^2 h^2)$ estimate on the area of the support of $\bar \partial F_\lambda$, we can get 
\[ \| \dot{f_\lambda (P)} \|_{L^2( \Gamma_{r,g} )} \leq C \lambda^{-1} h^{-1},\]
with choosing $N=2$.

Similarly the operator $\ddot{f_\lambda(P)}$ has an operator norm bound of $C M^{-4} h^{-4/3}$.

\end{proof}

Then we can state our main result in this section,

\begin{prop}\label{exterior}
One has
\begin{equation}\label{exterior estimate}
\left<f_\lambda(P)u_h, u_h \right>_{ \Gamma } \leq C \lambda^{-1}  \|u_h\|^2_{L^2(\Gamma) }.
\end{equation}
\end{prop}\

\begin{rem}
The main idea of proof is inspired by  \cite[Proposition 5.2]{BHT}. Since we work on the eigenfunctions in forbidden regions, which have exponentially small upper bound in $h$ rather than the ones in allowed regions, which have polynomial upper bound in $h$, each quantity in our proof should be treated differently.

Also note that the above estimate is different from \cite[Proposition 5.2]{BHT} since the right hand side quantity is restricted on a hypersurface.
\end{rem}

\begin{proof}
Let
\[L(r) = h^2 \left<f_\lambda(P)u_h, u_h \right>_{ \Gamma_{r,g} }. \]
With dots indicating differentiation with respect to $r$, and notice $\partial_r=\partial_{x_n}$
\[ \dot L( r)=h^2 \left<\dot {f_\lambda(P)} u_h, u_h \right>_{ \Gamma_{r,g} }+2h^2 \left<{f_\lambda(P)} u_h, \partial_{x_n} u_h \right>_{ \Gamma_{r,g} }  \]


\noindent and
\begin{align} \label{ex second}
\ddot L(r )=& h^2 \left<\ddot {f_\lambda(P)} u_h, u_h \right>_{ \Gamma_{r,g} }+ 4h \left<\dot{f_\lambda(P)} u_h, h\partial_{x_n} u_h \right>_{ \Gamma_{r,g}}  \nonumber \\
&+ 2 \left< f_\lambda(P) h\partial_{x_n}u_h, h \partial_{x_n} u_h \right>_{ \Gamma_{r,g} } + 2 \left< f_\lambda(P) u_h, h^2\partial^2_{x_n} u_h \right>_{ \Gamma_{r,g} }.
\end{align}

By Taylor expansion we write
\begin{equation}\label{Taylor}
L( r)=L(0) + r \dot L(0) +\int_0^r ds \int_0^s  \ddot L(t) dt.
\end{equation}

The key observation is that the term $ \left< f_\lambda(P) u_h, h^2\partial^2_{x_n} u_h \right>_{ \Gamma_{r,g} }$ in \eqref{ex second} is strongly positive, since 
\begin{align}\label{identity P}
h^2\partial^2_{x_n} u_h=  Pu_h &+h^2\frac{\partial_{x_n} |\det g|} {2|\det g|} \partial_{x_n} u_h+ ( V-E(h) )u_h
\end{align}
and 
\begin{equation}\label{positivity P}
P\geq \lambda h
\end{equation} on the support of $f_\lambda$. We will show that, unless $L(0)$ is very small, this term drives the polynomial growth in $h$ of $L$ which will contradict the value by integrating $L$ on the interval $[0, \lambda h]$.

Using Lemma \ref{lemma}, we estimate
\begin{equation}\label{d L}
\dot L(0)\geq - C \lambda^{-1} h \|u_h\|^2_{L^2( \Gamma )}.
\end{equation}

\noindent Similarly, for the terms on the right hand side of \eqref{ex second}, we have
\[ \vline\, \left< \ddot{f_\lambda(P)}u_h, u_h \right>_{ \Gamma_{r,g} } \vline \leq C \lambda^{-2} h^{-2} \|u_h\|^2_{L^2( \Gamma_{r,g} )},  \]

\begin{align*}
 4h \left<\dot{f_\lambda(P)} u_h, h \partial_{x_n} u_h \right>_{ \Gamma_{r,g} } 
 &= 2h^2 \frac{d}{dx_n} \left<\dot{f_\lambda(P)} u_h, u_h \right>_{ \Gamma_{r,g} } - 2h^2  \left<\ddot{f_\lambda(P)} u_h, u_h \right>_{ \Gamma_{r,g} } \\
 &\geq 2h^2 \frac{d}{dx_n} \left<\dot{f_\lambda(P)} u_h, u_h \right>_{ \Gamma_{r,g} } - C \lambda^{-2} \|u_h\|^2_{L^2( \Gamma_{r,g} )},
 \end{align*}
and induce the key inequality from \eqref{identity P} and \eqref{positivity P}
\begin{align}
 &2 \left< f_\lambda(P) u_h, h^2\partial^2_{x_n}u_h \right>_{ \Gamma_{r,g} } \nonumber \\
 =& 2 \Big< f_\lambda(P) u_h,  Pu_h +h^2\frac{\partial_{x_n} |\det g| } {2|\det g|} \partial_{x_n} u_h  +(V-E(h) )u_h \Big>_{ \Gamma_{r,g }  } \nonumber \\
 \geq& 2T h^{-2} L(r)+  \left< f_\lambda(P) u_h, h^2\frac{\partial_{x_n} |\det g|} {|\det g|} \partial_{x_n} u_h \right>_{ \Gamma_{r,g} } \nonumber \\
 \geq& 2T h^{-2} L(r)- Ch^2 \|u_h\|_{L^2( \Gamma_{r,g}) } \|\partial_{x_n} u_h\|_{L^2( \Gamma_{r,g}) }
 \end{align}
 in the second inequality notice that $P+V-E(h)\geq T>0$ on the support of $f_\lambda$, here $T$ is a small constant independent of $h$ and $r$. Also we applied Cauchy-Schwartz inequality in the last inequality.

From Lemma \ref{lem1}, for any fixed small $r$, taking $h$ sufficiently small, one has $\|u_h\|^2_{L^2( \Gamma_{r,g} )}< \|u_h\|^2_{L^2( \Gamma)}$. Setting $\tilde L(r)=\|u_h\|^2_{L^2( \Gamma_{r,g} )} $,  one also has $\frac{d}{d r} \tilde L(0)=2\left<u_h, \dot u_h\right>_ \Gamma =0$.  Hence $\frac{d}{d r} \tilde L(r)<0$ for any small $r$, which is equivalent to say that for any small $r$,
\begin{equation}\label{L2 compare}
\|u_h\|^2_{L^2( \Gamma_{r,g} )}< \|u_h\|^2_{L^2( \Gamma)}.
\end{equation}
With Lemma \ref{lem2} and the boundary condition $\|\partial_{x_n} u_h\|_{L^2(\Gamma)}=0 $, one can deduce that for any small $r$,
\begin{equation}\label{normal d compare}
\|\partial_{x_n}u_h\|^2_{L^2( \Gamma_{r,g} )}< C h^{-1} \|u_h\|^2_{L^2( \Gamma)}.
\end{equation}

Putting these together, we obtain
\begin{equation}
\ddot L(r)\geq - C \lambda^{-2}  \|u_h\|^2_{L^2( \Gamma )} + 2h^2 \frac{d}{dx_n} \left<\dot{f_\lambda (P)} u_h, u_h \right>_{ \Gamma_{r,g} }+ 2T h^{-2} L(r).
\end{equation}

Using this in \eqref{Taylor}, we get an inequality
\begin{align*}
L(r)\geq & L(0) - r C \lambda^{-1} h \|u_h\|^2_{L^2( \Gamma )}  \\
&+ \int_0^rds \int_0^s dt \left( - C \lambda^{-2} \|u_h\|^2_{L^2( \Gamma)} + 2h^2 \frac{d}{dt} \left<\dot{f_\lambda(P)} u_h, u_h \right>_{ \Gamma_{t, g} }+ 2T h^{-2} L(t) \right).
\end{align*}
For $r\in[0, \lambda h]$, by \eqref{L2 compare} and Lemma \ref{lemma}, the first two terms in the big bracket can be absorbed by the $rC \lambda^{-1} h \|u_h\|^2_{L^2(\Gamma )}$ term if $h$ is small enough. We get
\[ L(r)\geq L(0) - r C \lambda^{-1} h \|u_h\|^2_{L^2( \Gamma )} + 2T \int_0^rds \int_0^s h^{-2} L(t)  dt.\]

Then one has following comparison: $ L(r)\geq Z(r)$ where $Z(r)$ satisfies the corresponding equality
\[Z(r)= L(0) - r C \lambda^{-1} h \|u_h\|^2_{L^2(\Gamma )} + T \int_0^rds \int_0^s h^{-2} Z(t) dt, \quad Z(0)=L(0). \]
This we can solve exactly: differentiating twice gives us
\begin{equation}
\ddot{Z(r)}= T h^{-2} Z(r), \quad Z(0)=L(0), \,\, \dot Z(0)=-  C \lambda^{-1} h \|u_h\|^2_{L^2(\Gamma )}. 
\end{equation}

\noindent The solution is
\begin{align*}
Z(r)&=L(0)cosh( {\sqrt T} h^{-1} r) - \frac{1}{\sqrt T}C \lambda^{-1}  h^{2} \|u_h\|^2_{L^2(\Gamma )} sinh({\sqrt T} h^{-1}r) \\
&\geq \left( L(0)- \frac{1}{\sqrt T} C \lambda^{-1} h^{2 } \|u_h\|^2_{L^2(\Gamma )} \right) e^{{\sqrt T} h^{-1}r}
\end{align*}


\noindent Now suppose, for a contradiction, that $L(0)$ was bigger than $ \frac{1}{\sqrt{T}} C \lambda^{-1} h^{2 }  \|u_h\|^2_{L^2(\Gamma )}$. This would tell us that
\begin{equation}
L(r)\geq Z(r) \geq C \lambda^{-1} h^{2 }  \|u_h\|^2_{L^2(\Gamma)} e^{ {\sqrt T} h^{-1}r} \quad r\in [0, \lambda h].
\end{equation}

Integrating this on $[0,  \lambda h]$ gives
\begin{equation}\label{lower}
\int_0^{ \lambda h} L(r)dr \geq  C  \frac{e^{\sqrt{T}\lambda}-1}{ \lambda} h^3 \|u_h\|^2_{L^2(\Gamma )}.
\end{equation}

On the other hand by the definition of $L(r)$, one has
\[\int_0^{ \lambda h} L(r)dr \leq C \lambda h^{3}   \|u_h\|^2_{L^2(\Gamma )}, \]
which contradicts with \eqref{lower} for sufficiently large $\lambda$.
We conclude that
\[ L(0)\leq  C \lambda^{-1} h^2 \|u_h\|^2_{L^2(\Gamma )},\]
proving this proposition.

\end{proof}

\medskip

\bibliography{reference}
\bibliographystyle{alpha}

\end{document}